\documentstyle[12pt]{article}
\newcommand{\const}{\mathop{\rm const}\limits}

 \newcommand{\mod}{\mathop{\rm mod}\limits}

\newcommand{\argmax}{\mathop{\rm argmax}\limits}

\begin{document}
% \large
\begin{center}

{\bf EXACT ASYMPTOTIC FOR THE TAIL OF MAXIMUM  OF  }\\

\vspace{3mm}

{\bf SMOOTH RANDOM FIELD DISTRIBUTION. }\\

\vspace{3mm}

{\sc Ostrovsky E.}\\

\vspace{3mm}

{\it Department of Mathematics and Statistics, Bar \ - \ Ilan University,
59200, Ramat Gan, Israel.}\\
e-mail: \ galo@list.ru \\

\vspace{3mm}

{\sc Abstract}\\

\end{center}

\vspace{3mm}

 We obtain in this paper using the saddle \ - \ point method
the expression for the exact asymptotic for the tail of
maximum of smooth (twice continuous differentiable)
 random field (process) distribution.

\vspace{3mm}


\vspace{3mm}


\vspace{3mm}

{\bf 1. Introduction. Notations. Statement of problem.} \par

 Let $ (\Omega, {\cal F}, {\bf P} ) $ be a  probability space with expectation
$ {\bf E } $ and variance $ {\bf Var}. $  Let also $ D $ be a
open bounded  convex domain with compact closure $ [D] $
in the space $ R^d $ with $ (d \ - \ 1) $  \ - \ dimensional boundary
$ \partial D = [D] \setminus D $ and let $ \xi = \xi(x) = \xi(x, \omega), \ \omega \in \Omega, \ x \in D $ be a twice continuous
differentiable  on the set $ [D] $
with probability one random process (field in the case
$ d \ge 2) $ with the values on the real line:

$$
\xi: D \times \Omega \to R^1.
$$

 Let us denote

$$
M = M(\omega) = \max_{x \in D} \xi(x),  \  T_M(u) = {\bf P}(M > u). \eqno(1.1)
$$

{\bf Our goal is the calculation of the exact asymptotic as $ u \to \infty $
for the  tail \ - \ function $ T_M(u) $ of maximum distribution in the terms of some finite \ - \ dimensional distributions of  the considered field }
 $ \xi(x). $ \par
 Recall that by definition  the asymptotical expression
$$
T_M(u) \sim Y(u), \ u \to \infty
$$

is said to be {\it exact},  iff

$$
\lim_{u \to \infty} T_M(u)/Y(u) = 1.
$$

 The symbol $ \asymp $ will denote as usually the weak  relation: we
write $ f(\cdot) \asymp g(\cdot), \ \lambda \in \Lambda $
for two numerical functions  $ f(\lambda), \ g(\lambda) $
defined on the arbitrary set $ \Lambda $ iff

$$
0 < \inf_{\lambda \in \Lambda} f(\lambda)/g(\lambda) \le
\sup_{\lambda \in \Lambda} f(\lambda)/g(\lambda) < \infty.
$$

 It is easy to see that the case of maximum absolute value
$ M_1 = \max |\xi(x)| $ or  $ M_- = \min \xi(x) $ and $ z \to \ - \ \infty $
may be considered analogously. \par

 The estimations of the probability $ T_M(u) $ as $ u \to \infty $ are used in
the probability theory \cite{Ledoux1}, \cite{Talagrand1},  \cite{Talagrand2},
\cite{Ostrovsky3}, theory of random fields \cite{Kozachenko1}, \cite{Ledoux1},
\cite{Talagrand3}, \cite{Ostrovsky1}, statistics \cite{Ostrovsky1}, \cite{Ostrovsky2}, theory of Monte \ - \ Carlo method \cite{Ostrovsky1}, \cite{Ostrovsky2}, reliability theory \cite{Ostrovsky1}, theory of approximation \cite{Ostrovsky4} etc.\\

 The exact asymptotic for the tail $ T_M(u) $ for the Gaussian fields $ \xi(x) $
was obtained in \cite{Piterbarg1}; see also \cite{Adler1}.  The upper and low
{\it bounds} for $ T_M(u) $ was obtained in many publications ( \cite{Fernique1},
\cite{Dudley1}, \cite{Pizier1}, \cite{Ostrovsky1},  \cite{Ostrovsky2},
\cite{{Ostrovsky6}, \cite{Talagrand1}, \cite{Talagrand2}, \cite{Talagrand4} }  etc. \par

{\bf Another denotations. } Let

$$
\eta_{i,j}(x) = \frac{\partial^2 \xi(x)}{\partial x_i \ \partial x_j }, \
i,j = 1,2,\ldots,d
$$
be  the {\it Hessian} of the random field $ \xi = \xi(x), $

$$
\zeta(x) = \zeta(x, \omega) = | \ \det(\eta) \ |^{1/2},  \ K(d) =
(2 \pi)^{- d/2},
$$

$$
I(\lambda) = I(\lambda, \omega)= \int_D \zeta(x) \ \exp(\lambda \ \xi(x)) \ dx,
\eqno(1.2)
$$
where  the "great" parameter $ \lambda, \ \lambda \to \infty $  belongs to the sector $ S = S(\epsilon) $ in the {\it complex} plane:

$$
S(\epsilon) = \{ \lambda: \ | \ \arg(\lambda) \ | \le \pi/2 - \epsilon  \},
$$

where $ \epsilon $ be a fixed number in the interval $ (0, 1/2). $\par
 Further, we denote for the values $ \lambda \in S(\epsilon): $

$$
G(\lambda) = K(d) \ {\bf E} I(\lambda) = K(d) \ \int_D {\bf E} [ \zeta(x) \
\exp( \lambda \ \xi(x) )] \ dx  \eqno(1.3)
$$
(We used Fubini theorem).\\

\vspace{3mm}

{\bf 2. Assumptions}.\\

\vspace{3mm}

{\bf 1.} We assume that the considered  random field $ \xi(x) $ is non \ - \
degenerate  in the following sense. For arbitrary finite set
of  pair \ - \ wise different  elements $ x_k, \ k = 1,2, \ldots, n+m:
k  \ne l \  \Rightarrow x_k \ne x_l $ of the set $ D: x_k \in D $ and for
all the values

$$
( \vec{i,j} ) =  \{ i_r, j_r \}, \ r = 1,2, \ldots, m
$$

the random vector $  \vec{\theta}_{( \vec{i,j}) } =  \vec{\theta}= $

$$
  \{ \xi(x_1), \xi(x_2), \ldots, \xi(x_n); \eta_{i_1,j_1}(x_{n+1}), \eta_{i_2,j_2}(x_{n+2}), \ldots, \eta_{i_m,j_m}(x_{n + m}) \}
$$
has a bounded density of distribution

$$
f_{\vec{\theta}}(y_1,y_2, \ldots, y_{n+m})  =
f_{\vec{\theta}}(y_1,y_2, \ldots, y_{n+m}; \ x_1,x_2, \ldots,x_n, x_{n+1},
x_{n + 2}, \ldots, x_{n + m})
$$
with respect to the usually Lebesgue measure:

$$
V_{ (\vec{i,j}) } (x_1, x_2, \ldots, x_n, x_{n+1}, x_{n+2}, \ldots, x_{n + m} )
\stackrel{def}{=}
$$

$$
\sup_{y_1, y_2, \ldots, y_n, y_{n+1}, y_{n+2}, \ldots, y_{n+m}} f_{\vec{\theta}}(y_1,y_2, \ldots, y_n, y_{n+1}, y_{n+2}, \ldots, y_{n+m}) < \infty.
\eqno(2.1)
$$

 The condition (2.1) guarantee, by virtue of theorem of Ilvisaker,  that the (random) point of maximum  of the field $ \xi = \xi(x): $

$$
x_0 = \argmax_{x \in [D]} \xi(x)
$$

there exists, is unique, belongs to the open set $ D $ and is non \ - \ degenerate:

$$
\det \left(\eta_{i,j}(x_0) \right) = \det \left( \frac{\partial^2 \xi(x)}{\partial x_i \ \partial x_j } (x_0) \ \right) \ne 0. \eqno(2.3)
$$
 The two last properties might be understood with probability one. \par

\vspace{2mm}

{\bf 2.} We will suppose again that the fields $ \xi(x) $ and its Hessian
$ \eta_{i,j}(x) $ satisfy the so \ - \ called uniform Kramer's  condition. In detail, we write the "expectation" decomposition

$$
\xi(x) = a(x) + \xi^o(x), \ {\bf E} \xi^o(x) = 0,
$$

such that $ a(x) = {\bf E}\xi(x);  \ a(\cdot) \in C^2([D]) $  and

$$
\eta_{i,j}(x) = \frac{\partial^2 a(x)}{\partial x_i \ \partial x_j } + \eta_{i,j}^o(x), \ {\bf E} \eta_{i,j}^o(x) = 0.
$$
 We assume that

$$
 \forall \lambda \in R^1  \ \exists \ \exp(\phi(\lambda))
\stackrel{def}{=} \sup_{x \in [D]} {\bf E} \exp \left( \lambda \xi^o(x)  \right) < \infty, \eqno(2.4)
$$

and analogously suppose
$$
\max_{i,j} \ \sup_{x \in [D]} {\bf E} \exp
\left( \lambda \eta_{i,j}^o(x)  \right) < \infty. \eqno(2.5)
$$

 The conditions (2.4) and (2.5) imply, e.g., that the integral $ I(\lambda) $ there exists for all values $ \lambda \in R^1. $  Moreover, we can introduce
the so \ - \ called $ B(\phi) $ space (see, for instance, \cite{ Ostrovsky1},
\cite{ Ostrovsky6} ) and define the {\it natural} semi \ - \ distance  on the
set $ D, \ z_1, z_2 \in D $ by the formula

$$
 d(z_1,z_2) =  d_{\phi}(z_1,z_2) = ||\xi^o(z_1) \ - \  \xi^o(z_2)||B(\phi).
$$

 We must recall briefly for readers convenience some facts about the theory
of $ B(\phi) $ spaces. \par

 Let $ \phi = \phi(\lambda), \lambda \in (-\lambda_0, \lambda_0), \ \lambda_0 = const \in (0, \infty] $ be some even strong convex which takes positive values for positive arguments  continuous function, such that $ \phi(\lambda) = 0
\ \Leftrightarrow \lambda = 0; $

$$
|\lambda| \le 1 \ \Rightarrow C_- \ \lambda^2 \le \phi(\lambda) \le C_+ \
\lambda^2;
$$
$ C_-,C_+ =  \const,  0 < C_- \le C_+ < \infty; $

$$
  \lim_{\lambda \to \lambda_0} \phi(\lambda)/\lambda = \infty.
$$

 We denote the set of all these function as $ \Phi; \ \Phi =
\{ \phi(\cdot) \}. $ \par
 We say that the {\it centered} random variable (r.v) $ \zeta =
\zeta(\omega) $ belongs to the space $ B(\phi), $ if there exists some non \ - \ negative constant $ \tau \ge 0 $ such that

$$
\forall \lambda \in (-\lambda_0, \lambda_0) \ \Rightarrow
{\bf E} \exp(\lambda \zeta) \le \exp[ \phi(\lambda \ \tau) ].
$$
 The minimal value $ \tau $ satisfying the last inequality for all values
$ \lambda \in (-\lambda_0, \lambda_0) $ is called a $ B(\phi) \ $ norm
of the variable $ \zeta, $ write $ ||\zeta||B(\phi) = $

$$
\inf \{ \tau, \ \tau > 0: \ \forall \lambda: \
|\lambda| < \lambda_0  \ \Rightarrow
 {\bf E}\exp(\lambda \zeta) \le \exp(\phi(\lambda \ \tau)) \}.
 $$
 Notice that in the considered in this article case $ \lambda_0 = \infty.$ \par

 This spaces are very convenient for the investigation of the r.v. having a
exponential decreasing
 tail of distribution, for instance, for investigation of the limit theorem,
the exponential bounds of distribution for sums of random variables,
non-asymptotical properties, problem of continuous of random fields,
study of Central Limit Theorem in the Banach space etc.\par

  The space $ B(\phi) $ with respect to the norm $ || \cdot ||B(\phi) $ and
ordinary operations is a Banach space which is isomorphic to the subspace
consisted on all the centered variables of Orlich's space $ (\Omega,F,{\bf P}), N(\cdot) $ with $ N \ - $ function

$$
N(u) = \exp(\phi^*(u)) - 1, \ \phi^*(u) = \sup_{\lambda} (\lambda u -
\phi(\lambda)).
$$
 The transform $ \phi \to \phi^* $ is called Young \ - \ Fenchel transform.
The proof of considered assertion used the properties of saddle \ - \ point method and theorem of Fenchel \ - \ Moraux:
$$
\phi^{**} = \phi.
$$

 The next facts about the $ B(\phi) $ spaces are proved in \cite{Ostrovsky1},
 p. 19 \ - \ 40:

$$
{\bf A.} \ \zeta \in B(\phi) \ {\bf iff} \ {\bf E } \zeta = 0, \ {\bf and} \ \exists C = const > 0,
$$

$$
T(|\zeta|, \ u) \le \exp(-\phi^*(u/C)), u \ge 0,
$$
where $ T(|\zeta|, u)$ denotes the {\it tail} of
distribution of the r.v. $ \zeta: $ \par

$$
T(|\zeta|, \ u) =  {\bf P}(| \zeta| > u), \ u \ge 0,
$$
and this estimation is in general case asymptotically exact. \par
 Henceforth  $ C, C_j $ will denote the non \ - \ essentially positive
finite "constructive" constants. \par
 More exactly, if $ \lambda_0 = \infty, $ then the following implication holds:

$$
\lim_{\lambda \to \infty} \phi^{-1}(\log {\bf E} \exp(\lambda \zeta))/\lambda =
K \in (0, \infty)
$$

if and only if

$$
\lim_{u \to \infty} (\phi^*)^{-1}( |\log T(\zeta,u)| )/u = 1/K.
$$
 Here and further $ f^{-1}(\cdot) $ denotes the inverse function to the
function $ f $ on the left \ - \ side half \ - \ line $ (C, \infty). $ \par

{\bf B.} The function $ \phi(\cdot) $ may be "constructive" introduced by
the formula
$$
\phi(\lambda) = \phi_0(\lambda) \stackrel{def}{=} \log \sup_{x \in D}
 {\bf E} \exp \left(\lambda \xi^o(x) \right),
$$
 if obviously the family of the centered r.v. $ \{ \xi^o(x), \ x \in D \} $ satisfies the {\it uniform } Kramer's condition:
$$
\exists \mu \in (0, \infty), \ \sup_{x \in D} T( |\xi^o(x)|, \ u) \le
\exp(-\mu \ u), \ u \ge 0.
$$
 In this case
 we will call the function $ \phi(\lambda) = \phi_0(\lambda) $
a {\it natural } function. \par

 {\bf C.} We define
$$
\psi(r) = \psi_{\phi}(r) = r/\phi^{-1}(r), \ r \ge 2.
$$

 Let us introduce a new norm (the so-called "moment norm")
on the set of r.v. defined in our probability space by the following way: the space $ G(\psi), $ or, in the other words, Grand Lebesgue Space (GLS) $ G(\psi) = G(\psi_{\phi} ) $ consist, by definition, on all the {\it centered}  r.v.
$ \{\zeta \} $ with finite norm

$$
||\zeta||G(\psi) \stackrel{def}{=} \sup_{r \ge 2} |\zeta|_r/\psi(r), \
|\zeta|_r \stackrel{def}{=} {\bf E}^{1/r} |\zeta|^r.
$$

 It is proved that the spaces $ B(\phi) $ and $ G(\psi) $ coincides:$ B(\phi) =
G(\psi) $ (set equality) and both
the norm $ ||\cdot||B(\phi) $ and $ ||\cdot|| $ are equivalent: $ \exists C_1 =
C_1(\phi), C_2 = C_2(\phi) = const \in (0,\infty), \ \forall \xi \in B(\phi) $

$$
||\zeta||G(\psi) \le C_1 \ ||\zeta||B(\phi) \le C_2 \ ||\zeta||G(\psi).
$$

{\bf D.} The definition of GLS $ G(\psi) $ spaces
is correct still for the non-centered random
variables $ \zeta.$ If for some non-zero r.v. $ \zeta \ $ we have
$ ||\zeta||G(\psi) < \infty, $ then for all positive values $ u $

$$
{\bf P}(|\zeta| > u) \le 2 \ \exp
\left( - \phi^*(u/(C_3 \ ||\zeta||G(\psi)))   \right).
$$

and conversely if a r.v. $ \zeta $ satisfies Kramer's condition, then
$ ||\zeta||G(\psi) < \infty. $ \par

 Without loss of generality we can and will suppose

$$
\sup_{x \in D} [ \ ||\xi^o(x) \ ||B(\phi)] = 1,
$$
(this condition is satisfied automatically in the case of natural choosing
of the function $ \phi: \ \phi(\lambda) = \phi_0(\lambda) \ ) $
and that the metric space $ (T,d) $ relatively the so \- \ called
{\it natural} distance (more exactly, semi \ - \ distance)

$$
d(z_1,z_2) \stackrel{def}{=} ||\xi^o(z_1) \ - \ \xi^o(z_2)|| B(\phi)
$$
is complete. \par
 For example, if $ \xi(x) $ is a {\it centered} Gaussian field:
$ {\bf E} \xi(x) = 0, \ x \in [D] $ and is {\it normed: }

$$
 \max_{x \in [D] } {\bf Var } [\xi(x)] = 1
$$
with covariation function

 $ W(z_1,z_2) = {\bf E} [ \xi(z_1) \ \xi(z_2)], $ then
$ \phi_0(\lambda) = 0.5 \ \lambda^2, \ \lambda \in R, $ and

$$
d(z_1,z_2) = d_{\phi_0}(z_1,z_2)= ||\xi(z_1) - \xi(z_2)||B(\phi_0) =
$$

$$
 \sqrt{ \bf{Var} [ \xi(z_1) - \xi(z_2) ]} =
\sqrt{ W(z_1,z_1) - 2 W(z_1,z_2) + W(z_2,z_2) }.
$$

{\bf E.} Let us introduce for any subset $ V, \ V \subset D $ the so-called
{\it entropy } $ H(V, d, \epsilon) = H(V, \epsilon) $ as a logarithm
(natural)
of a minimal quantity $ N(V,d, \epsilon) = N(V,\epsilon) = N $
of a balls $ S(V, t, \epsilon), \ t \in V: $
$$
S(V, t, \epsilon) \stackrel{def}{=} \{s, s \in V, \ d(s,t) \le \epsilon \},
$$
which cover the set $ V: $
$$
N = \min \{M: \exists \{t_i \}, i = 1,2,…, M, \ t_i \in V, \ V
\subset \cup_{i=1}^M S(V, t_i, \epsilon ) \},
$$
and we denote also
$$
H(V,d,\epsilon) = \log N; \ S(t_0,\epsilon) \stackrel{def}{=}
 S(T, t_0, \epsilon), \ H(d, \epsilon) \stackrel{def}{=} H(T,d,\epsilon).
$$
 It follows from Hausdorf's theorem that
$ \forall \epsilon > 0 \ \Rightarrow H(V,d,\epsilon)< \infty $ iff the
metric space $ (V, d) $ is precompact set, i.e. is the bounded set with
compact closure.\par

 It is known (see, for example, \cite{Ostrovsky1}, \cite{Ostrovsky6}) that
if the following series converges:

$$
\sum_{n=1}^{\infty} 2^{-n} H \left(D,d,2^{-n} \right) < \infty, \eqno(2.6)
$$
then  a (non \ - \ centered) r.v. $ \beta = \max_{x \in D} \xi(x) $
belongs to the space $ B^+(\phi): $
$$
T_{|\beta|}(u) \le 2 \exp \left(-\phi^*(u/C) \right), \ u \ge 1.
\eqno(2.7)
$$
 The condition (2.6) holds if for example the so \ - \ called metric dimension
of the set $ D $ relative the distance $ d = d_{\phi} $ is finite:
$$
\kappa \stackrel{def}{=} \overline{\lim}_{\epsilon \to 0+}
\frac{ H(D,d,\epsilon)}{ |\log \epsilon|} < \infty.
$$

  Henceforth we will suppose also the condition (2.6) (and following the conclusion (2.7)) is satisfied. \par
 Note that more modern result in the terms of "majorizing  measures" or
equally in the terms of "generic
chaining" see in \cite{Talagrand1}, \cite{Talagrand2},
\cite{Talagrand3},\cite{Talagrand4},  \cite{Ostrovsky6}.\par

\vspace{4mm}

{\bf 3. Main result.} \par

\vspace{4mm}

{\bf Theorem 1.} {\it We assert under formulated above conditions: as}
$ \lambda \to \infty $ {\it uniformly in } $ \lambda \in S(\epsilon) $

$$
{\bf E} \ e^{\lambda M} \sim  K(d) \
 \lambda^{d/2} \ G(\lambda). \eqno(3.1)
$$

{\bf Proof.} Let us consider the integral $ I(\lambda). $  Using the classical
saddle \ - \ point method (see, e.g., \cite{Fedorjuk1}, chapter 2, section 4),
 we obtain  that with probability one

$$
I(\lambda) \sim K(d) \ \lambda^{-d/2} \ e^{\lambda M}. \eqno(3.2)
$$

The passing to the limit as $ \lambda \to \infty, \ \lambda \in S(\epsilon) $
here and further
may be proved on the basis of equality (2.7) and theorem of dominated convergence.\par
 We get taking the expectation of equality (3.2):

$$
 G(\lambda) / K(d) \sim \lambda^{-d/2} \ {\bf E} e^{\lambda \ M}. \eqno(3.3)
$$
The equality (3.3) is equivalent to (3.1).\par

{\bf Corollary 1.}  As long as

$$
{\bf E} \exp(\lambda M) \sim
\lambda \int_{0}
^{\infty} \exp(\lambda z) \ T_M(z) \ dz,
$$

we conclude: $ \lambda \to \infty \ \Rightarrow $

$$
\int_{0}^{\infty} \exp(\lambda z) \ T_M(z) \ dz \sim  R(\lambda),
$$

where

$$
R(\lambda) \stackrel{def}{=} \frac{\lambda^{-1 + d/2}}
{( 2 \pi)^{d/2} } \ \int_D Q(\lambda,x) \ dx. \eqno(3.4)
$$
 It is evident that as $ \lambda \to \infty $

$$
\int_{0}^{\infty} \exp(\lambda z) \ T_M(z) \ dz \sim
\int_{- \infty}^{\infty} \exp(\lambda z) \ T_M(z) \ dz.
$$

\vspace{2mm}

{\bf 4. Examples.} \par

\vspace{2mm}

 It is possible to verify that for the smooth Gaussian fields $ \xi(x) $ the
asymptotical equality (3.4) coincides with the classical results belonging to
Piterbarg  \cite{Piterbarg1} and Adler \cite{Adler1}. \par
 Thus, we consider further only the non \ - \ Gaussian case. Namely, suppose that as $ \lambda \to \infty, \ \lambda \in S(\epsilon) $

$$
R(\lambda) \sim C(R) \ \lambda^{\alpha}  \ \exp(\lambda^q /q ) \eqno(4.1)
$$
for some constants $ \alpha, C(R), q; \ C(R) \in (0, \infty), q > 1,
 \ \alpha \in (-\infty, \infty).$ \par
 Introduce the conjugate power $ p = q/(q \ - \ 1) $ and a function
$ \phi_p(\lambda) $ as follows:

$$
\phi_p(\lambda) = \lambda^2, \ |\lambda| \le 1;
$$

$$
\phi_p(\lambda) = |\lambda|^p, \ |\lambda| > 1.
$$

We assume in addition to the condition (4.1) that

$$
\sup_{x \in [D]} ||\xi^o(x)|| B(\phi_p) < \infty   \eqno(4.2)
$$

and moreover that

$$
a(x) = {\bf E} \ \xi(x) \in C^2([D]),  \ \xi(\cdot) \in C^2([D])
( \mod \ {\bf P}),
$$

$$
\sum_{n=1}^{\infty} 2^{-n} H \left(D,d_p,2^{-n} \right) < \infty, \eqno(4.3)
$$
where

$$
d_p(z_1, z_2) = ||\xi^o(z_1) \ - \ \xi^o(z_2) ||B(\phi_p)
$$
is the natural semi \ - \ distance on the set $ [D]$ between the points
$ z_1,  z_2 $ from the set $ [D]. $

 It follows from the main result of  \cite{Kozachenko1}, \cite{Ostrovsky2} that

$$
|| \max_{x \in [D]} \ \xi(x) \ ||G(\psi_p) < \infty,  \ \psi_p =
\psi_p = \psi_{\phi_p}(\cdot),
$$
or equally

$$
T_M(z) \le \exp  \left( - (z/C)^p  \right), \ z \ge 0,
$$
as long as

$$
\phi_p^*(\lambda)  \asymp \phi_q(\lambda), \ \lambda \in (-\infty, \infty).
$$

 Taking into account the following asymptotical equality
(see \cite{Fedorjuk1}, chapter 2, section 2):

$$
\int_0^{\infty} y^{\gamma} \ \exp ( \lambda y \ - \ y^p/p) \ dy \sim
(2 \pi)^{1/2} \ \lambda^{\Delta} \ \exp(\lambda^q/q),  \eqno(4.4)
$$
where $ \gamma = \const, $

$$
\Delta = \frac{2 \gamma + 2 \ - \ p}{ 2(p \ - \ 1) },
$$
we conclude that under considered conditions and
using Tauberian \ - \ Richter theorems

$$
{\bf P} (M > u) \sim (2 \pi)^{ \ - \ 1/2} \ C(R) \ u^{\alpha(p \ - \ 1) \ -
\ 1 + p/2} \ \exp( \ - \ u^p/p), \ \eqno(4.5)
$$
as $ u \to \infty. $ \par

%  \newpage

{99}

\end{document}